\documentclass{pcam_au}
\usepackage{color}

\begin{document}

\title{The Tricomi Equation}

\author{Gui-Qiang G. Chen}
\affiliation{Mathematical Institute, University of Oxford}

\maketitle

The \emph{Tricomi equation} is a second-order partial differential equation of
mixed elliptic-hyperbolic type for $u(x,y)$
with the form:
\begin{displaymath}
u_{xx} + xu_{yy}=0.
\end{displaymath}

It was first analyzed in the work by Francesco Giacomo Tricomi (1923)
on the well-posedness of a boundary value problem.
The equation is hyperbolic in the half plane $x<0$, elliptic in the half plane $x>0$,
and degenerates on the line $x=0$.
Its characteristic equation is
\begin{displaymath}
dy^2+x dx^2 =0,
\end{displaymath}
whose solutions are
\begin{displaymath}
y\pm \frac{2}{3}(-x)^{\frac{3}{2}}=C
\end{displaymath}
for any constant $C$, which are real for $x<0$.
The characteristics comprise two families of semicubical parabolas lying in the half plane $x<0$,
with cusps on the
line $x=0$.
This is of hyperbolic degeneracy, for which the two characteristic families coincide, {\it perpendicularly} to the
line $x=0$.

\medskip
For $\pm x>0$, set $\tau=\frac{2}{3}(\pm x)^{\frac{3}{2}}$. Then the Tricomi equation becomes
the classical elliptic or hyperbolic {\it Euler-Poisson-Darboux equation}:
\begin{displaymath}
u_{\tau\tau}\pm u_{yy}+\frac{\beta}{\tau}u_\tau=0.
\end{displaymath}
The index $\beta=\frac{1}{3}$ determines the singularity of solutions
near $\tau=0$, equivalently, $x=0$.


\medskip
Many important problems in fluid mechanics and differential geometry
can be reduced to corresponding problems
for the Tricomi equation, particularly \textit{transonic flow problems}
and \textit{isometric embedding problems}.
The Tricomi equation is a prototype of the \emph{generalized Tricomi equation}:
\begin{displaymath}
u_{xx} + K(x)u_{yy}=0.
\end{displaymath}
For a steady-state transonic flow in $\mathbb{R}^2$, $u(x,y)$ is the stream function
of the flow,
$K(x)$ and $x$ are functions of the velocity,
which are positive at subsonic and negative
at supersonic speeds, and $y$ is the angle of inclination of the velocity.
The solutions $u(x,y)$ also serve as entropy generators for entropy pairs of the potential
flow system for the velocity.
For the isometric embedding problem of two-dimensional Riemannian manifolds
into $\mathbb{R}^3$, the function $K(x)$ has the same sign as the Gaussian
curvature.

\medskip
A closely related partial differential equation is the \emph{Keldysh equation}:
\begin{displaymath}
xu_{xx} + u_{yy}=0.
\end{displaymath}

It is hyperbolic when $x<0$, elliptic when $x>0$,
and degenerates on the line $x=0$.
Its characteristics are
\begin{displaymath}
y\pm \frac{1}{2}(-x)^{\frac{1}{2}}=C
\end{displaymath}
for any constant $C$, which are real for $x<0$.
The two characteristic families are (quadratic) parabolas lying in the half plane $x<0$
and coincide, {\it tangentially} to the degenerate line $x=0$,
which is of parabolic degeneracy.
For $\pm x>0$,  the Keldysh equation becomes the
elliptic or hyperbolic {\it Euler-Poisson-Darboux equation}
with index $\beta=-\frac{1}{4}$, by
setting $\tau=\frac{1}{2}(\pm x)^{\frac{1}{2}}$.
Many important problems in continuum mechanics
can be also reduced to corresponding problems
for the Keldysh equation, particularly shock reflection-diffraction problems
in gas dynamics.

\end{document}